\documentstyle{nyj}

\input{diagrams}
\input{mssymb.sty}

\title{Lefschetz Fibrations with unbounded Euler Class} 
\author{Thilo Kuessner}  
\address{Mathematisches Institut, Universit\"at M\"unchen, Theresienstrasse 39, D-80333 M\"unchen, Germany}
 
\email{kuessner@mathematik.uni-muenchen.de , http://www.mathematik.uni-muenchen.de/$\sim$kuessner}



\newtheorem{thm}{Theorem}
\newtheorem{lem}{Lemma}
\newtheorem{cor} {Corollary}

\newtheorem{lam}{Lemma}

\newenvironment{proof}{{\it Proof:}\quad}{\hfill$\Box$}

\begin{document}

\begin{abstract}
We investigate the bounded cohomology of Lefschetz fibrations: we show that the
Euler class
of a genuine Lefschetz fibration with distinct vanishing cycles (of fiber genus
$\ge2$) is not bounded.
As a consequence, we exclude the existence of negatively curved metrics on 
Lefschetz fibrations with more than one singular fiber.

\end{abstract}
\maketitle
\tableofcontents

%%%% **** The text of the paper starts here **** %%%%

The bounded cohomology $H_b^*\left(X;{\Bbb Z}\right)$ is an invariant of topological spaces,
which was introduced by Gromov in his work about the 
simplicial volume and has since then shown to be useful also in group theory and dynamics of group actions.\\
A cohomology class $\beta\in H^*\left(X;{\Bbb Z}\right)$
is said to be bounded
if it is in the image of the natural map $H_b^*\left(X;{\Bbb Z}\right)\rightarrow H^*\left(X;{\Bbb Z}\right)$.
Among other results, Gromov proved that (real) characteristic classes in $H^*\left(BG^\delta;{\Bbb R}\right)$ are bounded, if $G^\delta$ is
an algebraic subgroup of $GL\left(n,R\right)$ equipped with the discrete topology.
This generalized the classical Milnor-Sullivan theorem which states
that Euler classes of flat affine bundles are bounded.\\
In this article, we consider the Euler class of Lefschetz fibrations.
A well-known theorem of Morita says that
the Euler class of a surface bundle, with fiber of genus $\ge2$,
is bounded. We prove a converse to Morita's theorem.
\begin{thm} A Lefschetz fibration, with regular fiber of genus $\ge2$,
has bounded Euler class if and only if all singular fibers have the same vanishing cycle.\end{thm}
As an application, we can exclude the existence of negatively curved metrics on a large number of Lefschetz fibrations.\\
\\
{\bf Corollary 2}: {\em If a Lefschetz fibration, with regular fiber of genus $\ge2$, admits a Riemannian metric with negative 
sectional curvature everywhere, then it has at most one singular fiber.}\\
\\
We recall that any finitely presented group $\Gamma$ can be realised as the fundamental group
of a Lefschetz fibration (\cite{do},\cite{bk}). If $\Gamma$ happens to be word-hyperbolic, we will actually show that $\pi_2 M_\Gamma\not =0$.
 
\section{Preliminaries} 
{\bf Lefschetz fibrations.} A smooth map $\pi:M\rightarrow B$ from a smooth (closed, oriented, connected) $4$-manifold $M$ to a smooth (closed, oriented, oriented)
2-manifold $B$ is said to be a Lefschetz fibration, if 
it is surjective and $d\pi$ is surjective except at finitely many critical points $\left\{p_1,\ldots,p_k\right\}=:C\subset M$, having the property that there are complex coordinate charts (agreeing with the orientations of $M$ and $B$), $U_i$ around $p_i$ and $V_i$ around $\pi\left(p_i\right)$, 
such that in these charts $f$ is of the form $f\left(z_1,z_2\right)=z_1^2+z_2^2$, see \cite{gs}. 
After a small homotopy the {\bf critical points are in distinct fibers}, we assume this to hold for the rest of the paper.

The preimages of
points in $B-\pi\left(C\right)$ are called regular fibers. It follows from the definition that all regular fibers are diffeomorphic and that the restriction $\pi^\prime:=\pi\mid_{M^\prime}:M^\prime\rightarrow B^\prime$ to {\bf $M^\prime:=\pi^{-1}\pi\left(M-C\right)$ is a smooth fiber bundle over $B^\prime:=B-\pi\left(C\right)$}.
 
Let $\Sigma_g$ be the regular fiber, a closed surface of genus $g$, and let, for an arbitrary point $*\in\Sigma_g$, be $Map_{g,*}$ the group of diffeomorphisms $f:\Sigma_g\rightarrow\Sigma_g$ with $f\left(*\right)=*$, modulo homotopies fixing $*$.
It is well-known, cf. \cite{mo}, that for any surface bundle one gets
a  monodromy $\rho:\pi_1M^\prime\rightarrow Map_{g,*}$, which fits into the commutative diagram
$$\begin{diagram}
&\ \ 1\longrightarrow\ \ \ \ &\pi_1\Sigma_g\ \ \ \ &\longrightarrow\ \ \ \pi_1M^\prime\ \ \longrightarrow&\ \ \ \ \pi_1B^\prime\longrightarrow 1\\
&\ \ \ \ &\dTo^{id}&\dTo^\rho&\dTo\\
&\ \ 1\longrightarrow&\pi_1\Sigma_g&\longrightarrow Map_{g,*}\longrightarrow &Map_g\longrightarrow 1
\end{diagram}$$
It follows from the local structure of Lefschetz fibrations that, for a simple loop $c_i$ surrounding $\pi\left(p_i\right)$ in $B$, its image under the monodromy,
$\rho\left(c_i\right)$, is the Dehn twist at some closed 
curve $v_i\subset\Sigma_g$. $v_i$ is called the 'vanishing cycle'.

We point out the following fact: for $\gamma\in\pi_1\Sigma_g\subset\pi_1M^\prime$, the pointed
mapping class $\rho\left(\gamma\right)$ is a mapping which twists some loop representing $\gamma\in\pi_1\left(\Sigma_g,*\right)$ once along itself back and forth, such that it is homotopic (but not base-point preserving homotopic) to the identity. If $\Sigma_g$ carries a hyperbolic metric, then
there is a representative of $\rho\left(\gamma\right)$ which
can be lifted to a hyperbolic isometry with axis $\tilde{\gamma}\subset{\Bbb H}^2$, mapping $\tilde{*}$ to $\gamma\left(\tilde{*}\right)$, for any lift $\tilde{*}$ of $*$.
 
{\bf Euler class of Lefschetz fibrations.} 
For a topological space $X$, and a rank-2-vector bundle $\xi$ over $X$, 
one has an associated Euler class 
$e\left(\xi\right)\in H^2\left(X;{\Bbb Z}\right)$.

If $\pi:M\rightarrow B$ is a Lefschetz fibration, we may consider the tangent bundle of the fibers, $TF$,  except at points of $C$, where this is not well defined. We get a rank-2-vector bundle $L^\prime$ over $M-C$ with euler class 
$e^\prime:=e\left(TF\right)\in H^2\left(M-C;{\Bbb Z}\right)$.\\
By a standard application of the Mayer-Vietoris sequence, there is an isomorphism
$i^*:H^2\left(M;{\Bbb Z}\right)\rightarrow H^2\left(M-C;{\Bbb Z}\right)$ induced by the inclusion. Hence, $e:=\left(i^*\right)^{-1}e^\prime\in H^2\left(M;{\Bbb Z}\right)$ is well-defined.
In what follows we will denote $e$ as the Euler class of the Lefschetz fibration $\pi:M\rightarrow
B$.
It is actually true (but we will not need it) that there exists a rank-2-vector bundle $\xi$ over $M$ such that $\xi\mid_{M-C}\simeq TF$. It is the pull-back of the universal complex line bundle,
pulled back via the map $f:M\rightarrow CP^\infty$ corresponding to $e\in H^2\left(M;{\Bbb Z}\right)$
under the bijection $H^2\left(M;{\Bbb Z}\right)\simeq\left[M,CP^\infty\right]$.

{\bf ${\Bbb S}^1$-bundles associated to surface bundles.}
For any surface bundle $\pi^\prime:M^\prime\rightarrow B^\prime$ we may, after fixing a Riemannian metric, consider $UTF$, the unit tangent bundle 
of the fibers. We consider the case that the fiber has genus $g\ge 2$.
Then this ${\Bbb S}^1$-bundle is, according to \cite{mo}, 
equivalent to the flat $Homeo^+\left({\Bbb S}^1\right)$-bundle with 
monodromy $\partial_\infty\rho$, where $\partial_\infty: Map_{g,*}\rightarrow Homeo^+\left({\Bbb S}^1\right)$ is constructed as follows.

Recall that $\pi_1\Sigma_g$ is word-hyperbolic, since $g\ge2$.
For $f\in Map_{g,*}$ let $f_*: \pi_1\left(\Sigma_g,*\right)\rightarrow\pi_1\left(\Sigma_g,*\right)$ be the induced map of fundamental groups, and $\partial_\infty f_*$ the extension of $f_*$ to 
the Gromov boundary $\partial_\infty\pi_1\left(\Sigma_g,*\right)$. It is well-known that 
$\partial_\infty f_*$ is a homeomorphism and that there is a canonical homeomorphism 
$\partial_\infty\pi_1\left(\Sigma_g,*\right)\simeq {\Bbb S}^1$. \\
If $\gamma\in ker\left(Map_{g,*}\rightarrow Map_g\right)\simeq\pi_1\Sigma_g$,
then $\partial_\infty\rho\left(\gamma\right)\in PSL_2R\subset Homeo^+\left({\Bbb S}^1\right)$ is a hyperbolic map whose both fixed points are the endpoints of the lift of a representative of $\gamma$ passing through $*$.\\
One should be aware that the extension of $UTF$ to $M-C$ is not flat: a loop surrounding
a singular fiber is trivial in $\pi_1\left(M-C\right)$ but its monodromy is a Dehn twist, giving a nontrivial homeomorphism of ${\Bbb S}^1$.\\

{\bf Bounded Cohomology.} It will be important for us to distinguish between bounded cohomology with integer coefficients,
$H_b^2\left(X;{\Bbb Z}\right)$, and bounded cohomology with real coefficients, $H_b^2\left(X;{\Bbb R}\right)$.
We refer to \cite{iva} for definitions.
To avoid too complicated notation, we use the following convention: for $\beta\in H^*\left(X;{\Bbb Z}\right)$, we denote $\overline{\beta}\in H^*\left(X;{\Bbb R}\right)$ its image under the canonical
homomorphism $H^*\left(X;{\Bbb Z}\right)\rightarrow H^*\left(X;{\Bbb R}\right)$.
Also, we will not distinguish between $H_b^*\left(X;{\Bbb R}\right)$ and $H_b^*\left(\pi_1X;{\Bbb R}\right)$.

A cohomology class $\beta\in H^*\left(X;{\Bbb Z}\right)$ is called bounded if it belongs
to the image of the canonical homomorphism $H_b^*\left(X;{\Bbb Z}\right)\rightarrow H^*\left(X;{\Bbb Z}\right)$.\\
We will use the following two facts. (A) is proved in Bouarich's thesis, see \cite{bo}.
(B) is proved in \cite{ger}.
 
(A): If $1\rightarrow N\rightarrow \Gamma\rightarrow G\rightarrow 1$ is an exact sequence of groups, then there is an exact sequence
$$0\rightarrow H_b^2\left(G;{\Bbb R}\right)\rightarrow H_b^2\left(\Gamma;{\Bbb R}\right)\rightarrow H_b^2\left(N;{\Bbb R}\right)^G\rightarrow H_b^3\left(G;{\Bbb R}\right).$$
 
(B): For any group $\Gamma$, there is an exact sequence, natural with respect to group homomorphisms, $$H^1\left(\Gamma;{\Bbb R}/{\Bbb Z}\right)\rightarrow
H_b^2\left(\Gamma;{\Bbb Z}\right)\rightarrow H_b^2\left(\Gamma;{\Bbb R}\right).$$

{\bf Universal Euler class (\cite{gh}).} There is a class $\chi\in H^2\left(Homeo^+{\Bbb S}^1;{\Bbb Z}\right)$ such that, for any
representation $\rho:\pi_1M\rightarrow Homeo^+{\Bbb S}^1$ associated to a surface bundle with Euler class
$e$, one has
$\rho^*\chi=e$. 
By the explicit construction in \cite{mo} or \cite{gh}, $\chi$ is bounded.\\
For a discrete group $\Gamma$, $H_b^2\left(\Gamma;{\Bbb Z}\right)$ "classifies" actions of $\Gamma$ on ${\Bbb S}^1$ (see \cite{gh2},thm.6.6). In particular (\cite{gh},p.35), for $\rho:\Gamma\rightarrow Homeo^+\left({\Bbb S}^1\right)$, $ \rho^*\chi=0$ holds 
if and only if all $\rho\left(\gamma\right)$ with $\gamma\in\Gamma$ have a common fixed point on ${\Bbb S}^1$. (Note that the original statement in \cite{gh} is mistaken and would actually imply the existence of two common fixed points.)

\section {Proof of Theorem 1}

In this section we derive theorem 1 from the following two lemmata, which will be proved in subsections 2.1 and 2.2.
\begin{lam} Let $\pi:M\rightarrow B$ be a Lefschetz fibration with monodromy $\rho$ and let $e\in H^2\left(M;{\Bbb Z}\right)$ 
be a cohomology class whose restriction to $M^\prime$ has 
a preimage $e_b^\prime$ in $H_b^2\left(M^\prime;{\Bbb Z}\right)\simeq H_b^2\left(\pi_1M^\prime;{\Bbb Z}\right)$.
Let $N:=ker\left(\pi_1M^\prime\rightarrow\pi_1M\right)$ and $e_N$ the 
restriction of $e_b^\prime$ to $N$.\\
Then $\overline{e}$ is bounded if and only if $e_N\in ker\left(H_b^2\left(N;{\Bbb Z}\right)\rightarrow H_b^2\left(N;{\Bbb R}\right)\right)$.\end{lam}
\begin{lam} Let $\Gamma$ be a group, $\mathcal{A}$ a (possibly infinite) set of generators of $\Gamma$ and $\rho:\Gamma\rightarrow Map_{g,*}$ a 
representation such that

a) for all $a\in\mathcal{A}$ the rotation number of $\partial_\infty\rho\left(a\right)$ is zero,

b) there is no common fixed point on ${\Bbb S}^1$, that is, 

there is no $x\in {\Bbb S}^1$ with $\partial_\infty\rho\left(a\right)\left(x
\right)=x$ for all $a\in{\mathcal{A}}$.\\
Then the Euler class of $\rho$ does not belong to the kernel of the canonical homomorphism $H_b^2\left(\Gamma;{\Bbb Z}\right)\rightarrow H_b^2\left(\Gamma;{\Bbb R}\right)$.\end{lam}
\begin{proof}[Proof of Theorem 1] 

We have a commutative diagram $$\begin{diagram}
&&1&&1&&1&&\\
&&\dTo&&\dTo&&\dTo&&\\
1&\rTo& \Gamma&\rTo& N&\rTo& V&\rTo&1\\
&&\dTo&&\dTo&&\dTo&&\\
&\rTo&\pi_1F&\rTo&\pi_1M^\prime&\rTo& \pi_1B^\prime&\rTo&1\\
\end{diagram}$$
with all rows and columns being exact sequences.

A few remarks are in order about well-definedness of the involved homomorphisms.
The second line is the long exact homotopy sequences of the surface bundle $M^\prime\rightarrow B^\prime$.
Inclusion maps $ker\left(N\rightarrow V\right)$ to $ker\left(\pi_1M^\prime\rightarrow\pi_1B^\prime\right)$,
 hence $$\Gamma:=ker\left(N\rightarrow V\right)\subset\pi_1F.$$
Clearly, the projection maps $N$ to $ker\left(\pi_1B^\prime\rightarrow\pi_1B\right)=V$. 
Surjectivity of this
homomorphism does not follow from the commutative 
diagram, but is easy to see geometrically. Indeed, each 
simple loop $c_i$ surrounding a puncture can be lifted to an element $\hat{c_i}\in N$,
just working in coordinate charts. For $g\in\pi_1B$, we fix some
lift $\hat{g}\in\pi_1M$. Then $\hat{g}\hat{c_i}\hat{g}^{-1}$ is an element of $N$, projecting to
$gc_ig^{-1}$. Since $V$ is generated by elements of the form $gc_ig^{-1}$, we have
surjectivity.

To prove theorem 1, assume that $e$ is bounded. The restriction of $e$ to $M^\prime$ is the Euler class of the surface bundle $M^\prime\rightarrow B^\prime$ with fiber of genus $\ge 2$. As $B$ has genus $\ge 2$, $B^\prime$ is hyperbolic. So we have that the restriction of $e$ to $M^\prime$ is bounded.

Then
$e_N\in ker\left(H_b^2\left(N;{\Bbb Z}\right)\rightarrow H_b^2\left(N;{\Bbb R}\right)\right)$, according to lemma 1.
Letting $e_{fib}$ be the restriction of $e_N$ to $\Gamma:=\pi_1\Sigma_g\cap N\subset N$, we conclude in particular 
$e_{fib}\in ker\left(H_b^2\left(\Gamma;{\Bbb Z}\right)\rightarrow H_b^2\left(\Gamma;{\Bbb R}\right)\right)$. However, as explained in the preliminaries, 
if we are given a hyperbolic metric on $\Sigma_g$,
then,
for any $\gamma\in\pi_1\Sigma_g$, $\rho\left(\gamma\right)$ can be represented by a mapping which lifts to a hyperbolic isometry of ${\Bbb H}^2$. (The axis of the hyperbolic isometry projects to a loop
representing $\gamma$.)
This implies that $\partial_\infty\rho\left(\gamma\right)$ has exactly two fixed points. Hence 
$\partial_\infty\rho\left(\gamma\right)$
is of rotation number zero and we may apply lemma 2 to conclude that there must exist a common fixed point $x\in {\Bbb S}^1$ for all $\partial_\infty\rho\left(\gamma\right)$, $\gamma\in\Gamma$. 

Since $\Sigma_g$ is a compact hyperbolic surface, 
$\pi_1\Sigma_g\subset PSL_2{\Bbb R}\subset Homeo^+\left({\Bbb S}^1\right)$ consists only of 
hyperbolic isometries.
But if $\partial_\infty\rho\left(\gamma_1\right)$ and $\partial_\infty\rho\left(\gamma_2\right)$ were hyperbolic isometries which had exactly one fixed 
point in common, then $\partial_\infty\rho\left(\gamma_1\gamma_2^{-1}\right)$ 
would be parabolic, contradicting the cocompactness of $\pi_1\Sigma_g$ in $PSL_2{\Bbb R}$.
Therefore all 
$\partial_\infty\rho\left(\gamma\right), \gamma\in\Gamma$, are boundary maps of hyperbolic isometries 
with the same pairs of fixed points. These hyperbolic isometries belong to a 1-parameter group. Hence,
by discreteness, the image of $\Gamma$ must be cyclic.

Since $\Gamma\subset\pi_1\Sigma_g$ corresponds to the vanishing cycles, this means that there is only one vanishing cycle.
\end{proof}

We close this section by proving some corollaries.
\begin{cor}: Let $\Gamma$ be a word-hyperbolic group and $M_\Gamma$ a Lefschetz fibration (with regular fiber of genus $g\ge2$)
with $\pi_1M_\Gamma=\Gamma$, which has at least two singular fibers. 
Then $\pi_2M_\Gamma\not=0$.\end{cor}
Remark: To any finitely presented group $\Gamma$, there exists some Lefschetz fibration $M_\Gamma$ with $\pi_1M_\Gamma=\Gamma$ (\cite{do},\cite{bk}).

\begin{proof} 
If $\pi_1M$ is word-hyperbolic, then 
$H^2_b\left(\pi_1M,{\Bbb Z}\right)\rightarrow H^2\left(\pi_1M,{\Bbb Z}\right)$ is surjective, by the Gromov-Mineyev 
theorem (\cite{gm}).

Assume $\pi_2M=0$.
Then, by the Hopf-identity, $$H^2\left(\pi_1M;{\Bbb Z}\right)
\simeq H^2\left(M;{\Bbb Z}\right)
/\pi_2M=H^2\left(M,{\Bbb Z}\right).$$ Thus, surjectivity of 
$H^2_b\left(\pi_1M,{\Bbb Z}\right)\rightarrow H^2\left(\pi_1M,{\Bbb Z}\right)$
would imply surjectivity of 
$H^2_b\left(M,{\Bbb Z}\right)\rightarrow H^2\left(M,{\Bbb Z}\right)$. In particular, the Euler class would be bounded.

By theorem 1, we conclude that all singular fibers would have the same vanishing cycle. If there were two singular fibers, the unique
vanishing cycle would bound two disks, pasting together to a nontrivial element of $\pi_2M$ after
connecting the vanishing cycles by some cylinder.
(To see that the union of the two disks is nontrivial in $\pi_2$, consider the 
local handlebody decomposition of $M$: it contains two 2-handles attached to 
tubular neighborhoods of
the vanishing cycles, i.e., to copies of
${\Bbb S}^1\times {\Bbb D}^2$, where ${\Bbb S}^1$ corresponds to the vanishing cycle. If there
where a 3-ball 
bounding the union of these two disks, it could be made transverse to ${\Bbb S}^1\times {\Bbb D}^2$, hence intersecting it in a surface bounded by the vanishing cycle. But the vanishing cycle is not 0-homologous in ${\Bbb S}^1\times {\Bbb D}^2$.)\end{proof}

\begin{cor}
If a Lefschetz fibration, with regular fiber of genus $\ge2$, 
admits a Riemannian metric with negative           
sectional curvature everywhere, then it has at most one singular fiber.\end{cor}
\begin{proof} This follows from corollary 1 because a Riemannian manifold $M$
with negative sectional curvature has word-hyperbolic fundamental group and moreover, from
the Cartan-Hadamard theorem, $\pi_iM=0$ for $i\ge2$.\end{proof}

\begin{cor}
A Lefschetz fibration with regular fiber of genus $g\le2$, or base of 
genus $h\le3$,
does not
admit a Riemannian metric with negative sectional curvature everywhere.\end{cor}
\begin{proof} 
Lefschetz fibrations over a base of genus $h$, with regular fiber of genus $g$ and exactly one singular fiber, exist if and only if there is some Dehn twist in $Map\left(\Sigma_g,*\right)$ which can be written as a product of $h$ commutators, which according to \cite{ko} is possible if and only if $g\ge 3$ and $h\ge2$.\end{proof}

\subsection{Criterium for bounded Euler class}

Recall that, for a Lefschetz fibration $\pi:M\rightarrow B$ with critical points $C$, 
$B^\prime:=B-\pi\left(C\right)$ and $M^\prime:=\pi^{-1}\left(B^\prime\right)$, we have a monodromy map 
$\rho:\pi_1M^\prime\rightarrow Homeo^+\left({\Bbb S}^1\right)$ 
with Euler class $e^\prime\in H_b^2\left(\pi_1M^\prime;{\Bbb Z}\right)$. We will consider the subgroup $N:=ker\left(\pi_1M^\prime\rightarrow\pi_1M\right)$ and will denote $e_N\in H_b^2\left(N;{\Bbb Z}\right)$ the Euler class of $\rho\mid_N$.

\begin{lem} Let $\pi:M\rightarrow B$ be a Lefschetz fibration with monodromy $\rho$ and let $e\in H^2\left(M;{\Bbb Z}\right)$ be a cohomology 
class whose restriction to $M^\prime$ has a preimage $e_b^\prime$ in $H_b^2\left(M^\prime;{\Bbb Z}\right)\simeq H_b^*\left(\pi_1M^\prime;{\Bbb Z}\right)$.
Let $N:=ker\left(\pi_1M^\prime\rightarrow\pi_1M\right)$ and $e_N$ the
restriction of $e_b^\prime$ to $N$.\\
Then $\overline{e}$ is bounded if and only if $e_N\in ker\left(H_b^2\left(N;{\Bbb Z}\right)\rightarrow H_b^2\left(N;{\Bbb R}\right)\right)$.\end{lem}

\begin{proof}
Let $i:M^\prime\rightarrow M$ be the inclusion. The pull-back $i^*\overline{e}$ is 
bounded, by assumption.

From boundedness of $i^*\overline{e}$ and the commutative diagram $$\begin{diagram}
H_b^2\left(M;{\Bbb R}\right) &\rTo^{i^*}&H_b^2\left(M^\prime;{\Bbb R}\right)\\
\dTo^{}& &\dTo^{}\\
H^2\left(M;{\Bbb R}\right)&\rTo^{i^*}&H^2\left(M^\prime;{\Bbb R}\right)\\
\end{diagram}$$ we see that $\overline{e}$ is bounded if and only if $\overline{e}_b^\prime\in H_b^2\left(M^\prime;{\Bbb R}\right)$ is in the image of \\
$i^*:H_b^2\left(M;{\Bbb R}\right)\rightarrow H_b^2\left(M^\prime;{\Bbb R}\right)$. \\
We consider the exact sequence $1\rightarrow N\rightarrow \pi_1M^\prime\rightarrow\pi_1M\rightarrow 1$, with $N:=ker i_*$. Bouarich's exact sequence (A) implies that $\overline{e}_b^\prime\in im\left(i^*\right)$ if and only if the restriction of $e_b^\prime$ to $N$ is trivial in the {\em bounded} cohomology of $N$.   \end{proof}

Using the exact sequence $1\rightarrow\Gamma\rightarrow N\rightarrow V\rightarrow 1$, 
stated at the beginning
of the proof of theorem 1, one actually can apply Bouarich's exact sequence 
$H_b^2\left(V;{\Bbb R}\right)\rightarrow H_b^2\left(N;{\Bbb R}\right)\rightarrow H_b^2\left(\Gamma;{\Bbb R}\right)\rightarrow 1$
and conclude 
that $e_N
\in ker\left(H_b^2\left(N;{\Bbb Z}\right)\rightarrow H_b^2
\left(N;{\Bbb R}\right)\right)$ if and only if $e_\Gamma\in ker\left(H_b^2\left(\Gamma;{\Bbb Z}\right)\rightarrow H_b^2
\left(\Gamma;{\Bbb R}\right)\right)$ and the preimage of $\overline{e_N}$ in $H_b^2\left(V;{\Bbb R}\right)$ (which then exists by Bouarich's sequence) is zero.

\subsection{Mapping class groups generated by maps of rotation number zero}

Let $\Sigma$ be a closed surface and $f:\Sigma\rightarrow\Sigma$ a homeomorphism. We denote $f_*:\pi_1\left(\Sigma,*\right)\rightarrow\pi_1\left(\Sigma,*\right)$ the induced homomorphism for a fixed base point $*$, and $\partial_\infty f_*:{\Bbb S}^1\rightarrow {\Bbb S}^1$ the induced homeomorphism of the 
Gromov-boundary $\partial_\infty\pi_1\Sigma\simeq {\Bbb S}^1$, as in chapter 1. Let $Fix\left(\partial_\infty f_*\right)=\left\{p\in {\Bbb S}^1:\partial_\infty f_*\left(p\right)\right\}$ be the set of fixed points 
on the Gromov-boundary.

\begin{lem}: 
Let $\Gamma$ be a group, $\mathcal{A}$ a (possibly infinite) set of generators of $\Gamma$ and $\rho:\Gamma\rightarrow Map_{g,*}$ a
representation such that 
                                                                    
a) for all $a\in\mathcal{A}$ the rotation number of $\partial_\infty\rho\left(a\right)$ is zero, and

b) there is no common fixed point on ${\Bbb S}^1$, that is, 

there is no $x\in {\Bbb S}^1$ with $\partial_\infty\rho\left(a\right)\left(x\right)=x$ for all $a\in{\mathcal{A}}$.\\
Then the Euler class of $\rho$ does not belong to the kernel of the canonical homomorphism $H_b^2\left(\Gamma;{\Bbb Z}\right)\rightarrow H_b^2\left(\Gamma;{\Bbb R}\right)$.
\end{lem}
\begin{proof} 
For $\gamma\in{\mathcal{A}}$ let $j_\gamma:{\Bbb Z}\rightarrow\Gamma$ be the homomorphism such that $j_\gamma\left(1\right)=\gamma$.
 By functoriality of the exact sequence (B) (section 1), we have a commutative diagram  
$$\begin{diagram}
\Pi_{\gamma\in {\mathcal{A}}}\mbox{H}^1\left({\Bbb Z};{\Bbb R}/{\Bbb Z}\right)&\rTo^{\simeq}&\Pi_{\gamma\in {\mathcal{A}}}\mbox{H}_b^2\left({\Bbb Z};{\Bbb Z}\right) &\rTo&\Pi_{\gamma\in {\mathcal{A}}}\mbox{H}_b^2\left({\Bbb Z};{\Bbb R}\right)\\
\uTo^{\Pi\ j_\gamma^*}& &\uTo^{\Pi\ j_\gamma^*}& &\uTo^{\Pi\ j_\gamma^*}\\
\mbox{H}^1\left(\Gamma;{\Bbb R}/{\Bbb Z}\right)&\rTo&\mbox{H}^2_b\left(\Gamma;{\Bbb Z}\right)&\rTo&\mbox{H}^2_b\left(\Gamma;{\Bbb R}\right),\\
\end{diagram}$$
where the isomorphism $$H_b^2\left({\Bbb Z};{\Bbb Z}\right)\simeq {\Bbb R}/{\Bbb Z} \simeq H^1\left({\Bbb Z}; {\Bbb R}/{\Bbb Z}\right)$$ follows from prop.\ 3.1.\ in \cite{gh}. 
%Hom\left(\Pi_{i\in I}Z\gamma_i,R/Z\right)\simeq \Pi_{i\in I}Hom\left(Z;R/Z\right)$$ lets us conclude that also $$\Pi_{i\in I}\ j_i^*:H^1\left(\Gamma;R/Z\right)\rightarrow\Pi_{i\in I}H^1\left(Z;R/Z\right)$$ is an isomorphism.

Let $e\in H_b^2\left(\Gamma;{\Bbb Z}\right)$ be the Euler class of $\rho$. Its image $j_\gamma^*e\in H_b^2\left({\Bbb Z};{\Bbb Z}\right)$ is the Euler class of the representation of ${\Bbb Z}$ mapping $1$ to $\rho\left(\gamma_i\right)$. 
By theorem A3 in \cite{gh}, $j_\gamma^*e$ is
mapped to the rotation number 
of $\rho\left(\gamma\right)$ under the isomorphism $H_b^2\left({\Bbb Z};{\Bbb Z}\right)\simeq {\Bbb R}/{\Bbb Z}$. The rotation number of $\rho\left(\gamma_i\right)$ is zero, hence $j_\gamma^*e=0$ for all $\gamma\in{\mathcal{A}}$.

Now assume that $e$, the Euler class of $\rho$, belonged to the kernel of the canonical homomorphism $H_b^2\left(\Gamma;{\Bbb Z}\right)\rightarrow H_b^2\left(\Gamma;{\Bbb R}\right)$. It follows that $e\in H_b^2\left(\Gamma;{\Bbb Z}\right)$ has a preimage
$$E\in H^1\left(\Gamma;{\Bbb R}/{\Bbb Z}\right).$$ 

Since $\mathcal{A}$ generates $\Gamma$, the homomorphism $\Pi_{\gamma\in {\mathcal{A}}}
j_\gamma^*:H^1\left(\Gamma;{\Bbb R}/{\Bbb Z}\right)\rightarrow \Pi_{\gamma\in {\mathcal{A}}}\mbox{H}^1\left({\Bbb Z};{\Bbb R}/{\Bbb Z}\right)$ is injective. Hence, $\Pi_{\gamma\in {\mathcal{A}}}
j_\gamma^*e=0$ implies $E=0$. Therefore, also $e=0$.

According to \cite{gh2},
this contradicts assumption b).\end{proof}

\subsection{Relation with simplicial volume}
For a closed, orientable manifold $M$,
we consider the simplicial volume $\parallel M\parallel$, defined in \cite{gro}.
It is well-known (\cite{gro}) that $\parallel M\parallel >0$ if and only if the fundamental class $\omega_M\in H^{dim\left(M\right)}\left(M;{\Bbb R}\right)$ is bounded.

\begin{lem}
Let $\pi:M\rightarrow B$ be a Lefschetz fibration with regular fiber $F$ of genus
$g\left(F\right)\not=1$. Let
$e\in H^2\left(M;{\Bbb Z}\right)$ be the Euler class. Then

$\parallel M\parallel >0$ if and only if
the cup-product $e\cup\pi^*\omega_B$ is bounded.\end{lem}
\begin{proof} Assume that $e\cup\pi^*\omega_B$ is bounded.
To show that the fundamental class is bounded, it suffices to show that $e\cup\pi^*\omega_B$ is a non-zero multiple of the fundamental class, i.e., that 
$$<e\cup\pi^*\omega_B,\left[M\right]>\not=0.$$
The proof of this inequality
is a minor generalisation of the argument in \cite{hk}.

We work with de Rham-cohomology.
Define $\pi_*:H^2\left(M\right)\rightarrow H^0\left(B\right)$ by 
$\pi_*=D_B^{-1}\pi_*D_M$, where $D_B$ resp.\ $D_M$ are the Poincare duality maps. One has \\
$<\pi^*\alpha\cup\beta,c>=<\alpha\cup\pi_*\beta,\pi_*c>$ for any $\alpha,\beta\in H^*\left(M\right), c\in H_*\left(M\right)$. 

$e\cup\pi^*\omega_B$ is a multiple of the volume form. Therefore its value on $\left[M\right]$ does
not depend
on the zero-volume set $\pi^{-1}\pi\left(C\right)$. Hence,
$$\mid e\cup \pi^*\omega_B\left(\left[M\right]\right)\mid=\mid \int_{M-\pi^{-1}\pi\left(C\right)}e\cup\pi^*\omega_B\mid
=\mid\int_{B-\pi\left(C\right)}\pi_*e\cup\omega_B\mid.$$
Using, for $b\in B$, $<\pi_*e,\left[ b\right]>=<\pi_*e,\pi_*\left[F\right]>=<e,\left[F\right]>=\chi\left(F\right)$,
we get $$\mid e\cup \pi^*\omega_B\left(\left[M\right]\right)\mid
=
\mid \chi\left(F\right)\int_B\omega_B\mid\not=0$$ because $\chi\left(F\right)\not=0$.\end{proof}

If genus$\left(B\right)\ge 2$, then $\pi^*\omega_B$ is bounded and, thus,
a sufficient 
condition for $\parallel M\parallel>0$ is boundedness of the Euler class $e$. One should not expect this condition to be necessary. To understand under which conditions boundedness of $e\cup\pi^*\omega_B$ (and hence nontriviality of $\parallel M\parallel$) holds, it would be necessary to understand the fourth bounded cohomology
of Lefschetz fibrations better. 

%The existing theory only admits partial results, a sample of them is the following.

%\begin{pro}
%Let $M$ be a Lefschetz fibration with regular fiber $F$ of genus $g\left(F\right)\ge2$, base $B={\Bbb S}^2$ and at least 3 singular fibers. Then $\parallel M\parallel > 0$.\end{pro}
%\begin{proof} The assumption of at least 3 singular fibers implies that $B^\prime$ has negative Euler characteristic.
%Hence the surface bundle $M^\prime\rightarrow B^\prime $ has positive simplicial volume and the restriction of $\overline{e}\cup\pi^*\omega_B$ to $M^\prime$ is bounded. 

\end{document}